\input amstex
\documentstyle{amsppt}
\magnification=\magstep1                        
\hsize6.5truein\vsize8.9truein                  
\NoRunningHeads
\loadeusm

\magnification=\magstep1                        
\hsize6.5truein\vsize8.9truein                  
\NoRunningHeads
\loadeusm

\document
\topmatter

\title
The number of real zeros of polynomials with constrained coefficients 
\endtitle

\author
Tam\'as Erd\'elyi
\endauthor

\address Department of Mathematics, Texas A\&M University,
College Station, Texas 77843 \endaddress

\email terdelyi\@tamu.edu
\endemail

\thanks {{\it 2020 Mathematics Subject Classifications.} 26C10, 11C08}
\endthanks

\dedicatory 
Dedicated to Edward Saff on the occasion of his 80th birthday
\enddedicatory 

\keywords location of zeros, polynomials with constrained coefficients
\endkeywords

\date September 30, 2024 \enddate

\abstract
We prove that there is an absolute constant $c > 0$ such that every 
polynomial $P$ of the form
$$P(z) = \sum_{j=0}^{n}{a_jz^j}\,, \quad |a_0| = 1\,, \quad |a_j| \leq M\,, \quad a_j \in \Bbb{C}\,, \quad M \geq 1\,,$$
has at most $cn^{1/2}(1+\log M)^{1/2}$ zeros in the interval $[-1,1]$. 
This result is sharp up to the multiplicative constant $c > 0$ and
extends earlier results of Borwein, Erd\'elyi, and K\'os from the case $M=1$ to the case $M \geq $1.
This has also been proved recently with the factor $(1+\log M)$ rather than $(1+\log M)^{1/2}$ 
in the Appendix of a recent paper by Jacob and Nazarov by using a different method.
We also prove that there is an absolute constant $c > 0$ such that every polynomial $P$ of the above form  
has at most $(c/a)(1+\log M)$ zeros in the interval $[-1+a,1-a]$ with $a \in (0,1)$.
Finally we correct a somewhat incorrect proof of an earlier result by Borwein and Erd\'elyi 
by proving that there is a constant $\eta  > 0$ such that every polynomial $P$ of 
the above form with $M = 1$ has at most $\eta n^{1/2}$ zeros inside any polygon 
with vertices on the unit circle, where the multiplicative constant $\eta > 0$ depends only on the polygon. 
\endabstract

\endtopmatter

\head 1. Introduction \endhead
We study the location of the zeros of polynomials of the form 
$$P(z) = \sum_{j=0}^{n}{a_jz^j}\,, \quad |a_0| = 1\,, \quad |a_j| \leq M\,, \quad a_j \in \Bbb{C}\,, \quad M \geq 1\,. \tag 1.1$$ 
For some of closely related results see the List of References and/or the Introduction in [BEK99].

\head 2. Notation \endhead

For a continuous function defined on an interval $I$ we define
$$\|f\|_{I} := \sup_{x \in I}{|f(x)|}\,.$$
For $z_0 \in \Bbb{C}$ and $r > 0$ let $D(z_0,r)$ denote the open disk of the complex plane
centered at $z_0$ of radius $r$. The closed disk centered at $z_0$ of radius $r$ will be
denoted by $\overline{D}(z_0,r)$.
For $M \geq 1$ let $\Cal{S}_M$ denote the collection of all analytic functions $f$ defined in the open unit disk $D(0,1)$
satisfying the growth condition
$$|f(z)| \leq \frac{M}{1-|z|}\,,\qquad z \in D(0,1)\,. \tag 2.1$$
In this paper $c_j$ denotes a suitable positive absolute constant.

\head 3. Results \endhead

Our first theorem gives an upper bound for the number of zeros of a polynomial $P$ of the form (1.1)
in the interval $[-1,1]$. This has been proved with the factor $(1+\log M)$ rather than 
$(1+\log M)^{1/2}$ in the Appendix of a recent paper [JN24] by Jacob and Nazarov, while 
the case $M=1$ was proved in [BEK99, Theorem 4.1]. We offer a proof in the case of $M \geq 1$ 
by using Jensen's formula (Lemma 4.1) and improving the dependence on $M$.

\proclaim{Theorem 3.1} 
There is an absolute constant $c > 0$ such that every polynomial $P$ of the form (1.1)
has at most $c(n(1+\log M))^{1/2}$ zeros in $[-1,1]$. 
\endproclaim

Our second theorem gives an upper bound for the number of zeros of a polynomial $P$ of the form (1.1) 
in the interval $[-1+a,1-a]$, $a \in (0,1)$. For $M=1$ the theorem below has been proved in 
[BEK, Theorem 4.2]. 

\proclaim{Theorem 3.2}
Let $a \in (0,1)$. There is an absolute constant $c > 0$ such that every polynomial $P$ of the form (1.1) 
has at most $(c/a)(1+\log M)$ zeros in $[-1+a,1-a]$.
\endproclaim

The sharpness of the above theorems may be compared with [BEK99, Theorem 2.3] stating that 
there exist polynomials $P$ of the form 
$$P(z) = \sum_{j=0}^{n}{a_jz^j}\,, \quad a_0 \neq 0, \quad |a_j| \leq 1\,, \quad a_j \in \Bbb{R}\,, \tag 3.1$$ 
such that $P$ has a zero at $1$ with multiplicity at least
$$\min \left\{\frac {1}{6} \left(n \, (1 - \log |a_0|)\right)^{1/2} - 1\,, n\right\}\,.$$
Multiplying by $M := 1/|a_0|$ in (3.1) shows that there exist polynomials $P$ 
of the form (1.1) such that $P$ has a zero at $1$ with multiplicity at least
$$\min \left\{\frac {1}{6} \left( n(1+\log M) \right)^{1/2} - 1\,, n\right\}\,.$$ 
It is also shown in [BEK99, Theorem 2.8] that if $|a_0| \geq \exp(-L^2)$ and $|a_j| \leq 1$ for each 
$j = L^2+1, L^2+2, \ldots, n$, then the polynomial
$$P(x) = \sum_{j=0}^n {a_j x^j}\,, \qquad  a_j \in {\Bbb C}\,,$$
has at most $\textstyle{\frac{44}{7}}(L+1)n^{1/2} + 5$ zeros at $1$.

Jacob and Nazarov [JN24] proved a general theorem recently including that there is an absolute 
constant $c>0$ such that there are polynomials of degree $n$ with coefficients in $\{-1,1\}$ 
having at least $cn^{1/2}$ {\it distinct} zeros in the interval $[-1,1]$; and there are polynomials of degree 
at most $n$ with coefficients in $\{0,1\}$ having at least $cn^{1/2}$ {\it distinct} zeros in the 
interval $[-1,1]$. 

The proof of our next theorem corrects a somewhat incorrect proof of an earlier result in [BE97, Theorem 2.1].

\proclaim{Theorem 3.3}
There is a constant $\eta > 0$ such that every polynomial of the form (1.1) with $M = 1$ has at most 
$\eta n^{1/2}$ zeros inside any polygon with vertices on the unit circle, where the multiplicative 
constant $\eta > 0$ depends only on the polygon.
\endproclaim

\head 4. Lemmas \endhead

One of our main tools is Jensen's formula stated as Lemma 4.1 below. See [BE95] for instance.

\proclaim{Lemma 4.1}
Suppose $h$ is a nonnegative integer, $z_0 \in \Bbb{C}$, and
$$f(z) = \sum_{j=h}^{m}{a_j(z-z_0)^j}, \quad a_j \in \Bbb{C}, \quad a_h \neq 0\,,$$
is analytic in an open disk centered at $z_0$ of radius greater than $R > 0$, and suppose that
the zeros of $f$ in $\overline{D}(z_0,R) \setminus \{z_0\}$ are $z_1,z_2, \ldots,z_m$,
where each zero is listed as many times as its multiplicity. Then
$$\log|a_h| + h\log R + \sum_{j=1}^{m}{\log\left( \frac{R}{|z_j|} \right)} = 
\frac{1}{2\pi} \int_0^{2\pi}{\log|f(Re^{i\theta})| \, d\theta}\,.$$
\endproclaim

Our second main tool is [BEK99, Theorem 5.1] stated below.

\proclaim{Lemma 4.2} There are absolute constants $c_1 > 0$ and $c_2 > 0$ such that
$$|f(0)|^{c_1/b} \leq \exp \left( \frac {c_2}{b} \right) \, \|f\|_{[1-b,1]}$$
for every $f \in \Cal{S}_1$ and $b \in (0,1]$.
\endproclaim

Lemma 4.2 implies the following by a simple linear scaling.

\proclaim{Lemma 4.3}
Let $I_a := [1-5a/4,1-a]$, $a \in (0,4/5]$. There is an absolute constant $c_3 > 0$ such that
$$(eM)^{-c_3/a} \leq \|P\|_{I_a}$$
for every polynomial $P$ of the form (1.1).
\endproclaim

\demo{Proof of Lemma 4.3} 
A simple estimation shows that every polynomial $P$ of the form (1.1) is contained in $\Cal{S}_M$.
It is also straightforward that $f \in \Cal{S}_M$ implies that $g_r$ defined by $g_r(z) := f(rz)$ 
with $0 < r \leq 1$ satisfies $g_r \in \Cal{S}_M$ and $g_r(0) = f(0)$. We apply this observation with
$r := 1-a$ and $b := a/4$ which implies that the interval $I_a = [1-5a/4,1-a]$ contains the interval 
$[(1-b)(1-a),1-a]$. Observe also that $f \in \Cal{S}_M$ implies that $h_M := f/M \in \Cal{S}_1$ and 
$h_M(0) = f(0)/M$ to conclude that the lemma follows from Lemma 4.2.
\qed \enddemo

To prove Theorem 3.1 we need the lemma below. 

\proclaim{Lemma 4.4} Assume that in addition to $M \geq 1$ we have $1+\log M \leq n$. 
Let $J_{n,M} := [1-n^{-1/2}(1+\log M)^{1/2},1]$. There is an absolute constant $c_4 > 0$ such that 
every polynomial $P$ of the form (1.1) has at most $c_4(n(1+\log M))^{1/2}$ zeros in the interval $J_{n,M}$.
\endproclaim

\demo{Proof of Lemma 4.4} 
Let $P$ be a polynomial of the form (1.1). It follows from Lemma 4.3 with $a := (4/5)n^{-1/2}(1+\log M)^{1/2}$ 
that there is a point $z_0 \in J_{n,M}$ such that
$$\log|P(z_0)| \geq -(5/4)c_3n^{1/2}(1+\log M)^{-1/2}\,(1+\log M) \geq  -2c_3n^{1/2}(1+\log M)^{1/2}\,. \tag 4.1$$
We will apply Jensen's formula (Lemma 4.1) with the disk $\overline{D}(z_0,2n^{-1/2}(1+\log M)^{1/2})$. 
Using the bounds for the coefficients of $P$ of the form (1.1), we get
$$|P(z)| \leq (n+1)M \big( 1+2n^{-1/2}(1+\log M)^{1/2} \big)^n \leq (n+1)M \exp(2n^{1/2}(1+\log M)^{1/2})$$  
for every $z \in \overline{D}(z_0,2n^{-1/2}(1+\log M)^{1/2})\,.$
Hence
$$\log|P(z)| \leq \log(n+1) + \log M + 2n^{1/2}(1+\log M)^{1/2} \tag 4.2 $$
for every $z \in \overline{D}(z_0,2n^{-1/2}(1+\log M)^{1/2})\,.$
Denote the zeros of $P$ in 
$$\overline{D}(z_0,n^{-1/2}(1+\log M)^{1/2})$$ 
by $z_1,z_2,\ldots,z_m$, where each zero is listed as many times as its multiplicity. 
Estimating by Jensen's formula (Lemma 4.1) with the disk $\overline{D}(z_0,2n^{-1/2}(1+\log M)^{1/2})$ 
and recalling (4.1) and (4.2), we obtain 
$$\split & -2c_3n^{1/2}(1+\log M)^{1/2} + 0 + m\log 2 \leq  
-2c_3n^{1/2}(1+\log M)^{1/2} + \sum_{j=1}^{m}{\log\left( \frac{2n^{-1/2}}{|z_j|} \right)} \cr 
\leq & \frac{1}{2\pi} \, 2\pi(\log(n+1) + \log M + 2n^{1/2}(1+\log M)^{1/2})\,. \cr \endsplit$$ 
Therefore
$$\split m & \leq \frac{1}{\log 2}\big( \log(n+1) + \log M + 2n^{1/2}(1+\log M)^{1/2} + 2c_3n^{1/2}(1+\log M)^{1/2} \big) \cr 
& \leq c_4n^{1/2}(1+\log M)^{1/2}\,, \cr
\endsplit $$
where the inequalities 
$$\log M = (\log M)^{1/2} (\log M)^{1/2} \leq n^{1/2}(\log M)^{1/2}$$
and $\log(n+1) \leq n^{1/2}$ have also been used.
Finally observe that the interval $J_{n,M}$ is contained in the disk 
$\overline{D}(z_0,n^{-1/2}(1+\log M)^{1/2})$, and the lemma follows. 
\qed \enddemo

The proof of our next lemma is quite similar to that of Lemma 4.4.

\proclaim{Lemma 4.5}  
Let $a \in (0,4/5]$ and $I_a := [1-5a/4,1-a]$. There is an absolute constant $c_5 > 0$ such that  
every polynomial $P$ of the form (1.1) has at most $(c_5/a)(1+\log M)$ zeros in the interval $I_a$.
\endproclaim

\demo{Proof of Lemma 4.5}
Let $P$ be a polynomial of the form (1.1). By Lemma 4.3 there is a point $z_0 \in I_a$ such that 
$$\log|P(z_0)| \geq -(c_3/a)(1+\log M)\,. \tag 4.3$$
We will apply Jensen's formula (Lemma 4.1) with the disk $\overline{D}(z_0,a/2)$. Observe that
$$|P(z)| \leq \frac{2M}{a}\,, \qquad z \in \overline{D}(z_0,a/2)\,,$$
hence
$$\log|P(z)| \leq \log(2M)+ \log(1/a)\,, \qquad z \in \overline{D}(z_0,a/2)\,. \tag 4.4$$
Denote the zeros of $P$ in $\overline{D}(z_0,a/4)$ by $z_1,z_2,\ldots,z_m$, where each zero is listed
as many times as its multiplicity. Estimating by Jensen's formula (Lemma 4.1) with the disk $\overline{D}(z_0,a/2)$ 
and recalling (4.3) and (4.4), we obtain
$$\split -(c_3/a)(1+\log M) + 0 + m\log 2 \leq & -(c_3/a)(1+\log M) + \sum_{j=1}^{m}{\log\left( \frac{a/2}{|z_j|} \right)} \cr 
\leq & \frac{1}{2\pi} \, 2\pi(\log(2M) + \log(1/a))\,, \cr \endsplit $$
hence
$$m \leq \frac{1}{\log 2}\big( \log(2M) + \log(1/a) + (c_3/a)(1+\log M) \big) \leq (c_5/a)(1+\log M)\,.$$
Finally observe that the interval $I_a = [1-5a/4,1-a]$ is contained in the disk $\overline{D}(z_0,a/4)$,
and the lemma follows. 
\qed \enddemo

\proclaim{Lemma 4.6}
Every polynomial $P$ of the form (1.1) has at most $\displaystyle{\frac{\log(2M)}{\log 2}}$ zeros in the disk $\overline{D}(0,1/4)$.
\endproclaim

\demo{Proof of Lemma 4.6}
Let $P$ be a polynomial of the form (1.1).
Let $z_1,z_2,\ldots,z_m$ be the zeros of $P$ in the disk $\overline{D}(0,1/4)$, where each zero is listed as many times as
its multiplicity. Using Jensen's formula (Lemma 4.1) with the disk $\overline{D}(0,1/2)$ we have
$$0 + 0 + m\log 2 \leq \sum_{j=1}^{m}{\log\left( \frac{1/2}{|z_j|} \right)} \leq \frac{1}{2\pi} \, 2\pi \log(2M)\,,$$
and the lemma follows.
\qed \enddemo

To prove Theorem 3.3 we need Lemmas 3.2, 3.3, and 3.5 from [BE97], which are formulated below as Lemmas 4.7, 4.8, and 4.9.
For the sake of completeness we present the proof of Lemmas 4.7, 4.8, and 4.9 below.9

\proclaim{Lemma 4.7}
Suppose $f \in \Cal{S}_1$, $\overline{D}(z_0,r) \subset \overline{D}(0,1)$ with $r \in [3/4,1)$, and 
the boundary circles have a common point. 
There is a constant $c(r) > 0$ depending only on $r \in [3/4,1)$ such that
$$\frac{1}{2\pi} \int_0^{2\pi}{\log |f(z_0 + re^{i\theta})| \, d\theta} \leq c(r)\,.$$
\endproclaim

\demo{Proof of Lemma 4.7}
This follows simply from Lemmas 4.10 and 4.11 and the growth condition (2.1) with $M=1$.
\qed \enddemo

The lemma below is a consequence of Jensen's formula (Lemma 4.1) and Lemma 4.2.  

\proclaim{Lemma 4.8}
There is an absolute constant $c_6 > 0$ such that every polynomial $P$ of the form (1.1) with $M=1$
has at most $c_6(n^{1/2} + nr)$ zeros in any disk $\overline{D}(u,r)$ with $u \in \Bbb{C}$, $|u|=1$, and $r > 0$. 
\endproclaim

\demo{Proof of Lemma 4.8}
Without loss of generality we may assume that $u=1$. We may also assume that  
and $n^{-1/2} \leq r \leq 1$ as the case $0 < r < n^{-1/2}$ follows from the case $r = n^{-1/2}$, 
and the case $r > 1$ is obvious.
Let $P$ be a polynomial of the form (1.1) with $M=1$. Observe that such a
polynomial satisfies the growth condition (2.1).
Choose a point $z_0 \in [1-r,1]$ such that
$$|P(z_0)| \geq \exp\left(\frac {-c_2}{r} \right)\,. \tag 4.5$$
There is such a point $z_0$ by Lemma 4.2.
Using the bounds for the coefficients of $P$ of the form (1.1), we have
$$\log |P(z)| \leq \log((n+1)(1+4r)^n) \leq \log(n+1) + 4nr\,, \qquad |z| \leq 1 + 4r\,. \tag 4.6$$
Let $m$ denote the number of zeros of $P$ in the disk $\overline{D}(z_0,2r)$.
Applying Jensen's formula (Lemma 4.1) with the disk $\overline{D}(z_0,4r)$ and recalling (4.5) and (4.6), we obtain
$$\frac {-c_2}{r} + m \log 2 \leq \log|P(z_0)| + m \log 2 \leq \frac{1}{2\pi}\,2\pi (\log(n+1) + 4nr)\,.$$
This, together with $n^{-1/2} \leq r \leq 1$, implies 
$$m \leq \frac{1}{\log 2} \left( \frac{c_2}{r} + \log(n+1) + 4nr \right) \leq c_6(n^{1/2} + nr)\,.$$ 
Now observe that $\overline{D}(1,r)$ is contained in $\overline{D}(z_0,2r)$, and the lemma follows.
\qed \enddemo

The lemma below is a consequence of Jensen's formula (Lemma 4.1) and Lemma 4.7. 

\proclaim{Lemma 4.9}
Suppose that $P$ is a polynomial of the form (1.1) with $M=1$. 
Let $\overline{D}(z_0,r) \subset \overline{D}(0,1)$, $r \in [3/4,1)$, $0 < \delta < r$, 
and the boundary circles have a common point. 
There is a constant $c(r) > 0$ depending only on $r \in [3/4,1)$ such that $P$ has at most $c(r)\delta^{-1}$ zeros in 
the disk $\overline{D}(z_0,r-\delta)$.
\endproclaim
 
\demo{Proof of Lemma 4.9} Let $P$ be a polynomial of the form (1.1) with $M=1$.
Observe that $r \in [3/4,1)$ and $\overline{D}(z_0,r) \subset \overline{D}(0,1)$ imply $|z_0| \leq 1/4$, hence
$$|P(z_0)| \geq 1-\sum_{j=1}^{\infty}{\left( \frac 14 \right)^j} = 1-1/3 = 2/3 \tag 4.7$$
Let $m$ denote the number of zeros of $P$ in the disk $\overline{D}(z_0,r-\delta)$.
Applying Jensen's formula (Lemma 4.1) with the disk $\overline{D}(z_0,r)$, (4.7), and Lemma 4.7, we obtain
$$\leqalignno{\log \left( \frac 23 \right) + m \log \left( \frac{r}{r-\delta} \right) 
& \leq \log |P(z_0)| + m \log \left( \frac{r}{r-\delta} \right) \cr
& \leq \frac{1}{2\pi}\int_0^{2\pi} {\log|P(z_0 + re^{i\theta})| \, d\theta} \leq c(r)\,. \cr}$$
As $3/4 \leq r < 1$ and $0 < \delta < r$, we have
$$\log \left( \frac{r}{r-\delta} \right) = \log \left( \frac{1}{1-(\delta/r)} \right) \geq \frac {\delta}{r} \geq \delta\,,$$
and, together with the previous inequality, this implies that $m \leq (c(r)+1)\delta^{-1}$.
\qed \enddemo

Lemmas 4.10 and 4.11 below are technical lemmas needed in the proof of Lemma 4.7. 
We study the closed disk centered at $(0,1)$ of radius $1$, and the closed disk centered at $(0,r)$ of radius $r \in (0,1)$  
in the Cartesian plane, so the larger disk contains the smaller disk, and the origin $(0,0)$ is the common point of 
the boundary circles.

\proclaim{Lemma 4.10}
Let $r \in (0,1)$. Let $d(x)$ be the distance between the points $(0,1)$ and $(x,r-(r^2-x^2)^{1/2})$, where $x \in [-r,r]$. 
Then $1-d(x) \geq c(r)x^2$ where $\displaystyle{c(r) = \frac{1-r}{2r}} > 0$ is a suitable constant depending only on 
$r \in (0,1)$. 
\endproclaim

\demo{Proof of Lemma 4.10}
Let $x \in [-r,r]$ and $r \in (0,1)$. We want to see that 
$$1-d(x) = 1 - \left( (x-0)^2 + ((1-r)+(r^2-x^2)^{1/2})^2 \right)^{1/2} \geq cx^2\,,$$
so it is sufficient to show that 
$$1-cx^2 \geq \left( x^2+(1-r)^2+r^2-x^2+2(1-r)(r^2-x^2)^{1/2} \right)^{1/2}\,,$$ 
so it is sufficient to show that
$$1-2cx^2 \geq (1-r)^2 + r^2 + 2(1-r)(r^2-x^2)^{1/2}\,,$$ 
so it is sufficient to show that 
$$cx^2 \leq (1-r)(r-(r^2-x^2)^{1/2}) = (1-r)\frac{x^2}{r + (r^2-x^2)^{1/2}}\,,$$
and this is true if $c = c(r) = (1-r)/(2r)$.
\qed \enddemo

\proclaim{Lemma 4.11}
Let $r \in (0,1)$. Let $d(x)$ be the distance between the points $(0,1)$ and $(x,r+(r^2-x^2)^{1/2})$, where $x \in [-r,r]$.
Then $1-d(x) \geq c(r)$, where $c(r) = r(1-r) > 0$ is a suitable constant depending only on $r \in (0,1)$.
\endproclaim

\demo{Proof of Lemma 4.11}
Let $x \in [-r,r]$ and $r \in (0,1)$. We have
$$\split 1-d(x) = & 1 - \left( (x-0)^2 + ((1-r)-(r^2-x^2)^{1/2})^2 \right)^{1/2} \cr
\geq & 1 - \left( 1-2r+2r^2-2(1-r)(r^2-x^2)^{1/2} \right)^{1/2} \cr 
\geq & 1 - ((1-2r+2r^2)^{1/2}) = \frac{1-1+2r-2r^2}{1+(1-2r+2r^2)^{1/2}} \cr
\geq & \frac{2r(1-r)}{1+1} = r(1-r)\,. \cr \endsplit$$ 
\qed \enddemo

\head 5.Proof of the Theorems \endhead

\demo{Proof of Theorem 3.2}
Let $P$ be a polynomial of the form (1.1). Without loss of generality we may assume that $a \in (0,1/4]$ 
otherwise the theorem follows from the case $a = 1/4$. 
By repeated applications of Lemma 4.5 we count the zeros of $P$ in the intervals
$$[1-(5/4)a,1-a]\,, \quad [1-(5/4)^2a,1-(5/4)a]\,, \quad \ldots\,, \quad [1-(5/4)^ka,1-(5/4)^{k-1}a]\,,$$
where $k$ is the smallest positive integer for which $(5/4)^ka \geq 3/4$. 
We obtain that $P$ has at most 
$$\sum_{j=0}^{k-1}{\frac{c_5(4/5)^j}{a} (1+\log M)} \leq \frac{5c_5}{a}(1+\log M)$$
zeros in the interval $[1/4,1]$. Combining this with Lemma 4.6 we obtain that $P$ has at most  
$$\frac{5c_5}{a}(1+\log M) + \frac{\log(2M)}{\log 2} \leq (c_7/a)(1+\log M)$$ 
zeros in the interval $[0,1-a]$. Replacing $P$ by $Q$ defined by $Q(z) = P(-z)$ we get that 
$P$ has at most $(c_7/a)(1+\log M)$ zeros in $[-1+a,0]$ and the theorem follows. 
\qed \enddemo

\demo{Proof of Theorem 3.1} Assume that in addition to $M \geq 1$ we have $\log M \leq n/4$, 
otherwise the theorem is trivial. The theorem follows from Lemma 4.4 and Theorem 3.2, as with the choice of 
$a := n^{-1/2}(1+\log M)^{1/2}$ we have 
$$(c/a)(1+\log M) = cn^{1/2}(1+\log M)^{-1/2}(1+\log M) = cn^{1/2}(1+\log M)^{1/2}\,.$$
\qed \enddemo

\demo{Proof of Theorem 3.3}
It is sufficient to prove the theorem for the number of zeros in a triangle with vertices $0$, $w$, and $w^{-1}$, 
where $w = e^{i\alpha}$ with $\alpha \in (0,\pi/3]$. 
Choose a 
$\gamma = \gamma(\alpha) \in (0,1/4]$ so that 
$$1-\gamma \geq \left| \gamma w - \left( \cos \alpha + \frac 12 (1-\cos \alpha) \right) \right|\,.$$    
Let
$$S_1 := \overline{D}(w,\beta n^{-1/2})\,, \qquad S_2 := \overline{D}(w^{-1},\beta n^{-1/2}),$$ 
$$S_3 := \overline{D}(\gamma w, 1-\gamma - n^{-1/2}), \qquad S_4 := \overline{D}(\gamma w^{-1}, 1-\gamma - n^{-1/2})\,.$$
Note that if $\beta = \beta(\alpha)$ and $n = n(\alpha)$ are sufficiently large then the triangle with 
vertices $0$, $w$, and $w^{-1}$ is covered by the union of $S_1$, $S_2$, $S_3$, and $S_4$. Hence the theorem 
follows from Lemmas 4.8 and 4.9. 
\qed \enddemo

\head 6. Acknowledgment \endhead
The author thanks Fedor Nazarov for the discussions on the topic and his calling my attention to his paper [JN24].

\enddocument